# Francis Guthrie's approach to the Four Color Problem

Asbjørn Brændeland

***Abstract***
Given a 5-critical graph *G*, for every edge $e = uv$ in *G* there is a 4-chromatic 5-subcritical graph $H = G - e$ with a color identical pair $\{u, v\}$. If *H* is non-planar so is *G*, but otherwise the only structure that can induce a color identity and keep things on the plane, is a color fixation chain of vertices and cycles, where *u* and *v* are separated by a cycle. So even if *H* is planar, the edge *e*, which joins *u* and *v* in *G*, must cross the cycle that separates them, making *G* non-planar.

Francis Guthrie put forth the Four Color Problem sometime around 1850: "Can the areas on any map be colored with at most four colors such that no pair of neighboring areas get the same color?"

Guthrie had drawn a map of four areas all bordering on each other, where four colors are required, but where also one area is enclosed by the three others and thus secluded from any added area, which could then be given the same color as the secluded one. We can call this *The Principle of the Seclusion of the Fourth Color* (allthough Francis probably did not use those words himself). From this he tried to prove that four colors actually do suffice, but as Francis' older brother Frederick notes 30 years later: "The proof which my brother gave did not seem altogether satisfactory to himself." Nevertheless, in spite of Francis' misgivings, we think he was on to something.

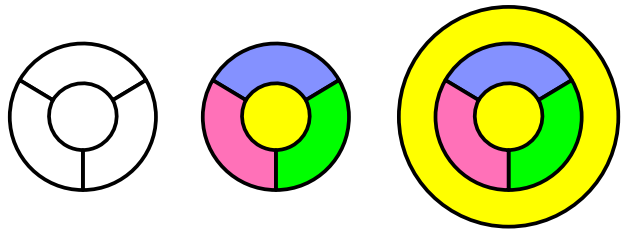

Figure 1

The dual of Francis' map is the *k*-complete graph, which is the smallest 4-critical graph.

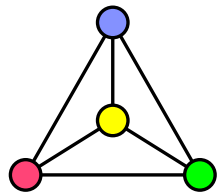

Figure 2

However, the set of 4-critical graphs has defied any attempt at characterization, so, trying to fit the principle of the seclusion of the fourth color into a characterization of the 4-critical graphs, simply won't work. To see where the principle fits in, we need to look at **subcritical** graphs and **color identical** pairs.

Consider a 4-critical graph *C* from which we remove an edge with end vertices *u* and *v*. The result *S* is a 4-subcritical, 3-chromatic graph. Now, in every 3-coloring of *S*, *u* and *v* must have equal colors, because if there were a 3-coloring of *S* in which *u* and *v* had different colors, that would also have been a 3-coloring of *C*, and *C* could not have been 4-critical. We say that *u* and *v* form a *color identical pair*.

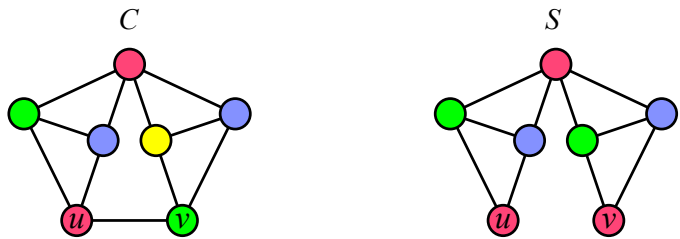

Figure 3

The subcritical graph *S* is a sequence of single vertices and 2-critical graphs (edges), where the coloring constraints are purely mechanical. When we color *u red*, the colors of the first edge are constrained to *green* and *blue*, the color of the middle vertex is constrained to *red*, the colors of the next edge to *green* and *blue* and the color of *v* to *red*. This gives the color identity {*u*, *v*} plus two others.

We say that a vertex constrained to a single color is ***fixed*** to that color and that a (*k* – 1)-critical subgraph of a *k*-chromatic graph constrained to a set of *k* – 1 colors is ***fixed*** to that color set, even if its vertices are not color fixed individually. *S* is then a ***color fixation chain***.

When, in a *k*-chromatic graph *G*, a single vertex *v* is adjacent to every vertex in a (*k* – 1)-critical subgraph *S* of *G*, *v* and *S* are mutually color fixed. Because of its adjacency to *v*, *S* cannot have more than *k* – 1 colors, and when *S* gets its color, *v* is fixed to the color *S* did not get, and when *v* gets its colors *S* is fixed to the remaining *k* – 1 colors.

It is possible for a vertex to be color fixed by an independent set of *k* – 1 vertices, provided the vertices in the set are fixed independently to *k* – 1 different colors. This works nicely with three colors, but with four colors one has to go above the plane. Notice that, also above the plane, and regardless of girth, the removal of an edge from a *k*-critical graph always creates a color identity. This is a totally different matter, though, since, without the mechanical coloring constraints of the triangles a color identity is given solely by the set of possible distributions of *k* colors among the cycles in the graph.

On the plane, however, the triangle rules. And here, the only way a single vertex can be color fixed in a 4-chromatic graph *G* is by being adjacent to every vertex in a 3-critical subgraph of *G*, i.e., an odd cycle.

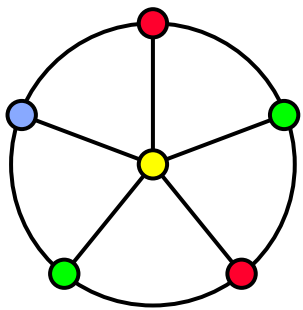

Figure 4

The vertex and the cycle above make out a wheel, and since the cycle is odd, we have an odd wheel, which is the simplest type of 4-critical graph—and also *the only type of 4-critical graph in which a single vertex can be color fixed*. This is very important, so we want to make it very clear.

For every 4-critical graph $G$, if $G$ has a vertex $v$ that is adjacent to every other vertex in $G$, then, removing $v$ leaves a face bounded by all the other vertices in $G$, and since $G - v$ is 3-chromatic, $G - v$ must be an odd cycle, and $G$ must be an odd wheel.

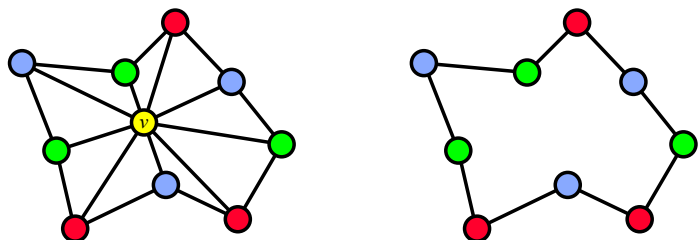

Figure 5

For every non-wheel 4-critical graph $G$ and for every vertex $v$ in $G$ there must be another vertex $u$, independent of $v$, in $G$ that can be given the same color as $v$. So there is *no* non-wheel 4-critical graph with a vertex that always gets a unique color, thus there is *no* non-wheel 4-critical graph with a vertex that can be color fixed.

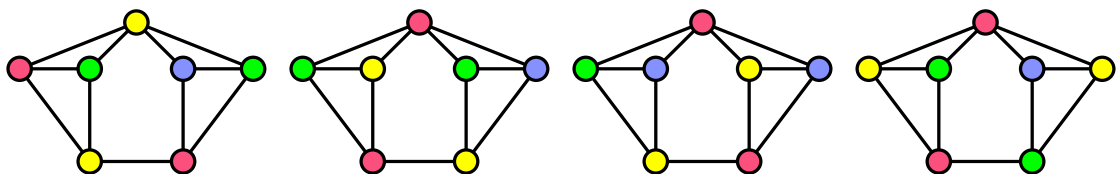

Figure 6

Given two vertices $u$ and $v$ that are both adjacent to every vertex in the same odd cycle. Regardless of how the cycle is colored, $u$ and $v$ must have equal colors, so $\{u, v\}$ is a color identical pair. But to avoid edge crossing and non-planarity, $u$ and $v$ must be on opposite sides of the cycle. This is the principle of the seclusion of the fourth color, applied to planar graphs.

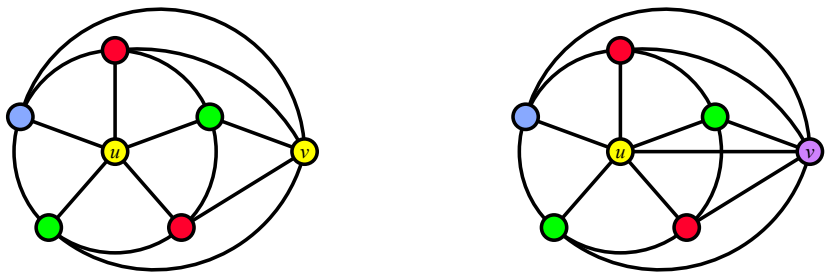

Figure 7

The two vertices and the cycle make out a pair of cycle-sharing wheels. Since the odd wheel is the only 4-critical graph where a single vertex can be color fixed, a pair of cycle sharing wheels is the only graph where exactly two vertices can be fixed to the same color. Thus a 5-critical graph with a subcritical that has exactly one color identical pair, must be a pair of cycle-sharing odd wheels with their hubs joined.

When two odd wheels share a hub their cycles are fixed to the same color set. Given such a pair, with the shared hub on the outside of both, place one new hub with spokes inside each cycle. This gives a color fixation chain of four wheels where the first share a cycle with the second, the second share a hub with the third and the third share a cycle with the fourth.

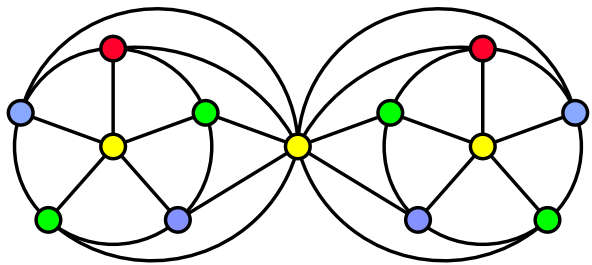

Figure 8

We can also describe this as a chain of alternating hubs and cycles, where all the hubs are fixed to the same color and all the cycles are fixed to the same color set. Every chain that starts and ends with a hub, is then the subcritical of a 5-critical graph.

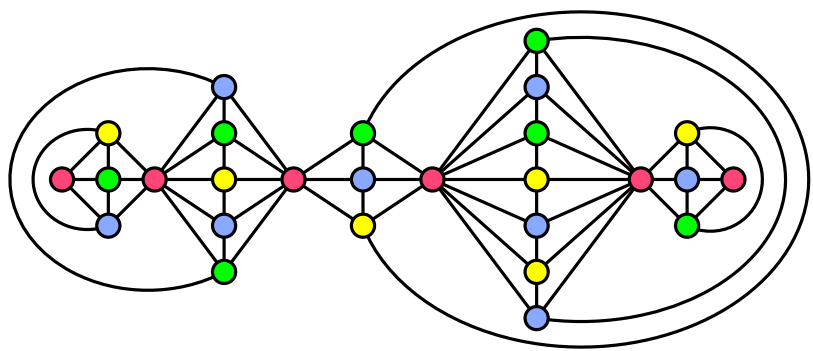

Figure 9

Every cycle in a 4-chromatic color-fixation chain divides the rest of the chain into two parts—the part inside and the part outside the cycle—thus, for every pair of color identical hubs in the chain, however close or distant, there is at least one cycle separating them.

So in the end, the diversity and unruliness of the 4-critical graphs has no bearing on the Four Color Problem, and Francis Guthrie had the right hunch after all.

The 5-subcritical color identity is itself contingent on the seclusion of the fourth color, which prevents the color identical vertices from being joined on the plane.